\documentclass[english]{amsart}
\usepackage{amssymb}
\usepackage{amsmath}
\usepackage{amscd}
\usepackage{amsthm}
\usepackage{amsfonts}
\usepackage{graphicx}
\usepackage[all]{xy}

\theoremstyle{plain}
\newtheorem{stam}{STAM}[subsection]
\newtheorem{thm}[stam]{Theorem}

\newtheorem{definition}[stam]{Definition}
\numberwithin{equation}{subsection} 

\theoremstyle{plain}
\newtheorem{conjecture}[stam]{Conjecture} 

\theoremstyle{plain}
\newtheorem{prop}[stam]{Proposition} 

\newtheorem{rmk}[stam]{Remark}

\DeclareMathOperator{\aff}{aff}
\DeclareMathOperator{\inter}{int}
\newcommand{\lfloorr}{\left\lfloor}
\newcommand{\rfloorr}{\right\rfloor}
\newcommand{\lbrackk}{\left\lbrack}
\newcommand{\rbrackk}{\right\rbrack}

\DeclareMathOperator{\conv}{conv}
\newcommand{\fR}{\mathbb{R}}
\DeclareMathOperator{\supp}{supp}

\newcommand{\calc}{\mathcal{C}}

\def\bbr{\mathbb{R}}
\def\suml{\sum\limits}

\makeatother
\begin{document}

\title{Some variations of Tverberg's Theorem}
\large

\author{Micha A. Perles, Moriah Sigron}

\begin{abstract}
 Define $T(d,r)=(d+1)(r-1)+1$. A well known theorem of Tverberg states that if $n\geq T(d,r)$, then one can partition any set of $n$ points in $\fR^d$ into $r$ pairwise disjoint subsets whose convex hulls have a common point. The numbers $T(d,r)$ are known as Tverberg numbers. Reay added another parameter $k$ $(2\leq k \leq r)$ and asked: what is the smallest number $n$, such that every set of $n$ points in $\fR^d$ admits an $r$-partition, in such a way that each $k$ of the convex hulls of the $r$ parts meet. Call this number $T(d,r,k)$. Reay conjectured that $T(d,r,k)=T(d,r)$ for all $d,r$ and $k$. In this paper we prove Reay's conjecture in the following cases: when $k\geq \lbrack \frac{d+3}{2}\rbrack$, or when $d<\frac{rk}{r-k}-1$, and for the specific values $d=3,r=4,k=2$ and $d=5,r=3,k=2$.

\end{abstract}

\maketitle
\newsavebox{\rb}
\savebox{\rb}(10,40)[bl]{
  \put(0,2.5){\line(0,1){40}}
  \put(0,0){\circle{5}}
  \put(-2.5,10){\line(1,0){5}}
  \put(3,7){$b$}
  \put(-2.5,35){\line(1,0){5}}
  \put(3,32){$r$}
}
\newsavebox{\br}
\savebox{\br}(10,40)[bl]{
  \put(0,2.5){\line(0,1){40}}
  \put(0,0){\circle{5}}
  \put(-2.5,10){\line(1,0){5}}
  \put(3,7){$r$}
  \put(-2.5,35){\line(1,0){5}}
  \put(3,32){$b$}
}
\newsavebox{\ry}
\savebox{\ry}(10,40)[bl]{
  \put(0,2.5){\line(0,1){40}}
  \put(0,0){\circle{5}}
  \put(-2.5,10){\line(1,0){5}}
  \put(3,7){$y$}
  \put(-2.5,35){\line(1,0){5}}
  \put(3,32){$r$}
}
\newsavebox{\by}
\savebox{\by}(10,40)[bl]{
  \put(0,2.5){\line(0,1){40}}
  \put(0,0){\circle{5}}
  \put(-2.5,10){\line(1,0){5}}
  \put(3,7){$y$}
  \put(-2.5,35){\line(1,0){5}}
  \put(3,32){$b$}
}
\newsavebox{\yb}
\savebox{\yb}(10,40)[bl]{
  \put(0,2.5){\line(0,1){40}}
  \put(0,0){\circle{5}}
  \put(-2.5,10){\line(1,0){5}}
  \put(3,7){$b$}
  \put(-2.5,35){\line(1,0){5}}
  \put(3,32){$y$}
}
\newsavebox{\yr}
\savebox{\yr}(10,40)[bl]{
  \put(0,2.5){\line(0,1){40}}
  \put(0,0){\circle{5}}
  \put(-2.5,10){\line(1,0){5}}
  \put(3,7){$r$}
  \put(-2.5,35){\line(1,0){5}}
  \put(3,32){$y$}
}
\newsavebox{\simplebar}
\savebox{\simplebar}(10,40)[bl]{
  \put(0,2.5){\line(0,1){40}}
  \put(0,0){\circle{5}}
}
\newsavebox{\rightmost}
\savebox{\rightmost}(10,40)[bl]{
  \put(0,2.5){\line(0,1){40}}
  \put(0,0){\circle{5}}
  \put(-2.5,35){\line(1,0){5}}
  \put(3,32){$p_{5,1}$}
}
\setcounter{tocdepth}{1}
\section{Introduction}

A well known theorem of Radon says that any set of $d+2$ or more points in $\mathbb{R}^{d}$ can be partitioned into two disjoint parts whose convex hulls meet. This follows easily from the fact that every set of d+2 points in $\mathbb{R}^{d}$ is affinely dependent.

The corresponding statement for partitions into more than two parts
is known as Tverberg's theorem.

\begin{thm}
\label{thm:Tverberg}
(H. Tverberg, 1966) Let $a_{1},\ldots,a_{n}$ be points in $\mathbb{R}^{d}$.
If $n >(d+1)(r-1)$ then the set N=\{1,\ldots,n\} of indices
can be partitioned into r disjoint parts $N_{1},\ldots,N_{r}$ in such
a way that the r convex hulls $\conv \{ a_{i}:i\in N_{j}\}$ (j=1,\ldots,r)
have a point in common.
\end{thm}
(This formulation covers also the case where the points $a_{1},\ldots,a_{n}$
are not all distinct.) Henceforth we use the abbreviation $a(N_{j})$
for $\{ a_{i}:i\in N_{j}\}$. The original proof (see \cite{T66}) was quite difficult.
In 1981 Tverberg published another proof, much simpler than the original one (see \cite{T81}).  Sarkaria \cite{S} gave a quite accessible proof, with some algebraic flavor. It seems that the simplest proof so far is due to Roudneff \cite{Ro}. See \cite{M} \S8.3 for further information.

The numbers $T(d,r)=(d+1)(r-1)+1$ are known as Tverberg numbers.  The condition $n\geq T(d,r)$ in Tverberg's theorem is extremely tight.  If $n<T(d,r)$, then almost always, for any $r$-partition $N_{1},\ldots,N_{r}$ of the set $N=\{1,\ldots,n\}$, even the intersection of the \textbf{affine} hulls $\aff(a(N_{j}))\,(j=1,\ldots,r)$ is empty.

In fact, there exists a polynomial $P$, not identically \underbar{}zero, in $n\cdot d$ scalar variables, $P(\vec{x_{1}},\ldots,\vec{x_{n}})=P(x_{11},\ldots,x_{1d},\ldots,x_{n1},\ldots,x_{nd})$, such that, for any $r$-partition $N_{1},\ldots,N_{r}$ of $N$, $\cap_{j=1}^{r}\aff(a(N_{j}))=\emptyset$ unless $P(a_{1},\ldots,a_{n})=0$. (For details, see \cite{PS}.)

In this paper we weaken the condition $\cap_{j=1}^{r}\conv a(N_{j})\neq\emptyset$
in Tverberg's theorem and ask only that each $k$ of
the convex hulls $\conv a(N_{j})\ (j=1,\ldots,r)$ meet, where $k$ is an additional
parameter, $2\leq k\leq r$. This weakened condition may perhaps require
fewer than $T(d,r)$ points. Thus we define $T(d,r,k)$ to be the
smallest positive integer $n$ with the following property: for any
list $a_{1},\ldots,a_{n}$ of points in $\mathbb{R}^{d}$ there is an
$r$-partition $N_{1},\ldots,N_{r}$ of the set of indices $N=\{1,\ldots,n\}$,
such that every $k$ of the $r$ convex hulls $\conv a(N_{j})$ have a
point in common.

The function $T(d,r,k)$ is clearly monotone non-decreasing in each
of the parameters $d,r,k$, and $T(d,r,r)=T(d,r)$.

If $r>d+1$, and each $d+1$ of the convex hulls $\conv a(N_{j})\,(j=1,\ldots,r)$
meet, then they all meet, by Helly's theorem. Thus $T(d,r,k)=T(d,r)$
for $d+1\leq k\leq r$. This reduces the interesting range of $k$
to $2\leq k\leq\min(r-1,d)$.

John R.  Reay (see \cite{Re}) settled the case $d=2$, showing that
$T(2,r,2)=T(2,r)$ for all $r\geq2$. He also showed that $T(3,3,2)=T(3,3)\,(=9)$
and made the following bold conjecture.

\begin{conjecture}
$T(d,r,k)=T(d,r)$ for all $2\leq k\leq r$.
\end{conjecture}

The meaning of Reay's conjecture is : If $n<T(d,r)$ then there exists a set $X \subset \bbr^d$, $|X|=n$, such that for every $r$-partition of $X$ there are two parts whose convex hulls are disjoint.

We don't really believe this is true. To press our point, consider
the case $d=r=1000$. By Tverberg's theorem, a million points in $\mathbb{R}^{1000}$
can be partitioned into one thousand parts whose convex hulls have
a common point. Is there a set of 999,999 points in $\mathbb{R}^{1000}$
that cannot be partitioned into 1000 parts whose convex hulls intersect
just pairwise? Seems implausible.

Nevertheless, the purpose of this paper is to establish parts of Reay's
conjecture. We show, by means of suitable examples, that Reay's conjecture
does hold in the following cases (Theorems \ref{thm:mainThm} - \ref{thm:532example}):

\begin{thm}
\label{thm:mainThm}
For every dimension $d \geq 2$ and for every
$r (\geq \lbrack \frac{d+3}{2} \rbrack)$:
\[ T(d,r, \lbrackk \frac{d+3}{2} \rbrackk)=T(d,r)=(d+1)(r-1)+1 .\]
\end{thm}

In particular, this shows that $T(3,4,3)=T(3,4)=13$. For $d=3,r=4,k=2$ we have the following:

\begin{thm}
\label{thm:342example}
$T(3,4,2)=T(3,4)=13$.
\end{thm}

Another class of cases is covered by:

\begin{thm}
\label{thm:secondThm}
For every $2\leq k<r$ and for every dimension $d<\frac{kr}{r-k}-1$:
\[ T(d,r,k)=T(d,r)=(d+1)(r-1)+1.\]
\end{thm}

Therefore, if $r=3$ and $k=2$ then $T(d,r,k)=T(d,r)$ provided $d<5$. The case $d=5$ is covered by the following:

\begin{thm}
\label{thm:532example}
$T(5,3,2)=T(5,3)=13$.
\end{thm}

in all cases, the examples are variations, specializations or perturbations of the following:
$d+1$ rays that emanate from the origin and positively span $\bbr^d$, with $r-1$ points chosen on each ray.

In order to put the ranges of Theorems ~\ref{thm:mainThm} and ~\ref{thm:secondThm} on the same scale, we can regard $k$ as the independent variable. Theorem ~\ref{thm:mainThm} establishes Reay's conjecture in the domain of $d$'s $d+1 \leq 2k-1$ (of which the subdomain $d+1 \leq k$ is trivial, in view of Helly's Theorem). For $d+1 \geq 2k$, Theorem ~\ref{thm:secondThm} establishes Reay's conjecture for $k<r< \frac{d+1}{d+1-k}k$. This domain of $r$'s reduces to $k<r<2k$ when $d+1=2k$, it shrinks with increasing $d$, and vanishes altogether for $d+1\geq k(k+1)$.

\section{Proof of theorem \ref{thm:mainThm}}
For the proof we will use the following (counter) example: let $p_{0},p_{1},\ldots,p_{d}\in\mathbb{R}^{d}$ be the vertices of a
$d$-simplex centered at the origin, i.e., $\sum_{i=0}^{d}p_{i}=\underbar{0}$
and each $d$ of the points $p_{0},p_{1},\ldots,p_{d}$ are linearly independent.
Let $D=\{0,1,\ldots,d\}$, and for $i\in D$ define $R_{i}=\{\lambda p_{i}:\lambda>0\}$ (the open ray emanating from \underbar{0} through $p_{i}$).

On each ray $R_{i}$ we choose $r-1$ distinct points.  The chosen points form a set $X\subset\mathbb{R}^{d}$, $|X|=(d+1)(r-1)=T(d,r)-1$.  We  show that in every partition of $X$ into $r$ parts ($X=C_{1}\cup \cdots\cup C_{r}$)
there is $j \leq \lfloorr \frac{d+3}{2} \rfloorr $ and there are parts $C_{i_1},\ldots,C_{i_j}$,
whose convex hulls have empty intersection. This will show that $T(d,r,k)=T(d,r)$
for $\lfloorr \frac{d+3}{2}\rfloorr \leq k\leq r$. We start with some preliminaries concerning the {}``positive basis'' $P=\{ p_{0},p_{1},\ldots,p_{d}\}$ of $\mathbb{R}^{d}$.

\subsection{Properties of the spanning set $P=\{ p_{0},p_{1},\ldots,p_{d}\}$}
\begin{prop}
Every point $x\in\mathbb{R}^{d}$ has a representation
\begin{equation}
\label{equ:positRep}
x=\xi_{0}p_{0}+\xi_{1}p_{1}+\cdots+\xi_{d}p_{d}
\end{equation}
 where $\min\{\xi_{0},\xi_{1},\ldots,\xi_{d}\}=0$. This representation
is \underbar{unique}.
\end{prop}
\begin{proof}
The vectors $p_{0},p_{1},\ldots,p_{d}$ span $\mathbb{R}^d$ linearly. In
fact, each $d$ of them form a linear basis of $\mathbb{R}^{d}$ .
Let $x=\sum_{i=0}^{d}\alpha_{i}p_{i}$ be some fixed representation
of $x$ in terms of $P$. The only linear dependences among
$p_{0},p_{1},\ldots,p_{d}$ are $\sum_{i=0}^{d}\lambda p_{i}=0,\,\lambda\in \mathbb{R}$.
Therefore the most general representation of $x$ in terms of $P$
is $x=\sum_{i=0}^{d}(\alpha_{i}-\lambda)p_{i}$, $\lambda\in \mathbb{R}$.
To obtain a representation with the smallest coefficient equal 0,
we must choose $\lambda=\min\{\alpha_{i}:i\in D\}$.
\end{proof}
We call (~\ref{equ:positRep}) the non-negative representation of $x$ (in terms
of $P$).
The \underbar{support} of $x$, (with respect to $P$) is defined
by \[
\supp x=\{ i\in D:\xi_{i}>0\}.\]
 Simple properties of $\supp x$:

\begin{enumerate}
\item $\emptyset\subseteq \supp x \subsetneq D$.
\item $\supp x=\emptyset$ iff $x=\underbar{0}$.
\item $\supp p_{i}=\{ i\}$.
\item $\supp \lambda x=\supp x$ for $\lambda>0$.
\item $\supp (x+y)\subseteq \supp x\cup \supp y$, with equality iff
$\supp x\cup \supp y\neq D$.
\item If $x\neq\underbar{0}$, then $\supp x\cup \supp(-x)=D$.
\end{enumerate}
Recall that our set $X$ consists of $r-1$ distinct points on each
ray $R_{i}$ $(i\in D)$. For a subset $C\subseteq X$, define $I(C)=\{ i\in D:C\cap R_{i}\neq\emptyset\}$.
Now make the following observations:

\begin{prop}
\label{prop:supp_x}
If $C\subseteq X$ and $x\in \conv C$, then $\supp x\subseteq I(C)$.
(This is obviously true also when I(C)=D.)
\end{prop}
When $I$ is a subset of $D$, we shall denote by $R(I)$ the union
$\cup\{ R_{i}:i\in I\}$.

\begin{prop}
\label{prop:CcapR(supp_x)}
Suppose $C\subseteq X$, and $x\in \conv C$. If $I(C)\neq D$ then $x\in \conv\{ C\cap R(\supp x)\}$.
\end{prop}
\begin{proof}
Suppose $x=\sum_{\nu=1}^{n}\gamma_{\nu}c_{\nu}$, where $c_{\nu}\in C,\,\gamma_{\nu}>0,\,\sum_{\nu=1}^{n}\gamma_{\nu}=1$.
If $c_{\nu}=\lambda_{\nu}p_{i}\in R_{i},\,\lambda_{\nu}>0$, then
$p_{i}$ will appear with a positive coefficient in the non-negative
representation of $x$ in terms of $P$, and therefore $i\in \supp x$.
Note that we have used the fact that $I(C)\neq D$.
\end{proof}
For points $a=\alpha p_{i}\in R_{i}$, $b=\beta p_{i}\in R_{i},\,(\alpha,\beta>0)$
we say that $a$ is \textbf{\underbar{lower}} than $b$ (or $b$ is \textbf{\underbar{higher}} than $a$) on $R_{i}$
if $\alpha<\beta$ (or, equivalently, if $\left\| a \right\| < \left\| b \right\|$).

\begin{prop}
\label{prop:ClowerThanC'}
Suppose $I\subsetneq D$. Let $C,C'$ be two finite subsets of
$R(I)(=\cup\{ R_{i}:i\in I\}).$
If, for each $i\in I$, every point of $C\cap R_{i}$ is lower (on
$R_{i}$) than every point of $C'\cap R_{i}$, then $\conv C \cap \conv C'=\emptyset$.
\end{prop}
\begin{proof}
Assume, w.l.o.g., that \textbar{}$I$\textbar{}=$d$. (We do \underbar{not}
assume that $C\cap R_{i}\neq\emptyset$ and $C'\cap R_{i}\neq\emptyset$
for all $i\in I$.) For each $i\in I$ choose a point $s_{i}=\sigma_{i}p_{i}\in R_{i}$
that is higher (on $R_{i}$) than every point of $C\cap R_{i}$ and
lower than every point of $C'\cap R_{i}$. The $d$ points $s_{i}$
($i\in I$) are linearly independent, and their affine hull $H=\aff\{ s_{i}:i\in I\}\subset\mathbb{R}^{d}$
is a hyperplane that does not pass through the origin. Denote by $H_{-},H_{+}$
the two open half spaces determined by H, and assume \underbar{0} $\in H_{-}$.
From our assumptions it follows that $C\subset H_{-}$ and $C'\subset H_{+}$,
hence $\conv C\cap \conv C'=\emptyset$.
\end{proof}

\begin{prop}
\label{prop:C1..Cn}
Suppose $U\subsetneq D$. Let $C_{1},C_{2},\ldots,C_{n}$ ($n \geq 2$)
be subsets of $X$. Assume
\begin{enumerate}
\item $\cap_{\nu=1}^{n}I(C_{\nu})\subseteq U$.
\item $I(C_{\nu})\subsetneq D$ for $\nu=1,2$.
\item For each $i\in U$, each point of $C_{1}\cap R_{i}$ is lower (on
$R_{i}$) than every point of $C_{2}\cap R_{i}$.
\end{enumerate}

Then $\cap_{\nu=1}^{n}\conv C_{\nu}=\emptyset$.
\end{prop}
\begin{proof}
Assume, on the contrary, that
$\cap_{\nu=1}^{n}\conv C_{\nu}\neq\emptyset$,
and suppose that
$x\in\cap_{\nu=1}^{n}\conv C_{\nu}.$
By Proposition~\ref{prop:supp_x} we conclude that
$\supp x\subseteq\cap_{\nu=1}^{n}I(C_{\nu})\subseteq U$.
Applying Proposition~\ref{prop:CcapR(supp_x)} to $C_{1}$ and
$C_{2}$, we
find that
\begin{equation}
\label{equ:x_in_conv}
x\in \conv(C_{\nu}\cap R(U))\,\, \text{for} \,\nu=1,2.
\end{equation}
Now invoke proposition~\ref{prop:ClowerThanC'} with $C=C_{1}\cap R(U)$,
$C'=C_{2}\cap R(U)$,
and $I=U$, to conclude that $\conv(C_{1}\cap R(U))\cap \conv(C_{2}\cap R(U))=\emptyset$,
which contradicts ($\ref{equ:x_in_conv}$).\end{proof}

\subsection{Completion of the proof of Theorem \ref{thm:mainThm}}

Let $X \subset \mathbb{R}^d$ be the set described at the beginning of this section
($r-1$ points on each of the rays $R_0,R_1,\ldots,R_d$),
and let $C_1,\ldots,C_r$ be an arbitrary partition of X into $r$ disjoint sets. Our aim is to
apply Proposition~\ref{prop:C1..Cn} to some $n$ of the parts $C_i$,
with $n$ as small as possible. We shall be able to do this
with some $n\leq \lbrackk \frac{d+3}{2} \rbrackk$.

Assume the parts $C_i$ are ordered in such a way that
\begin{enumerate}
    \item $|C_1| \leq |C_i|$ for $i=2,3,\ldots,r$.
    \item $|C_2 \cap R(I(C_1))| \leq |C_i \cap R(I(C_1))|$ for $i=3,4,\ldots,r$.
\end{enumerate}

From condition (1) we have
 \[ |I(C_1)| \leq |C_1| \leq \lbrackk \frac{1}{r}|X| \rbrackk = \lbrackk \frac{r-1}{r}(d+1) \rbrackk ,\]
 and therefore $|I(C_1)| \leq |C_1| \leq d$.

Condition (2) yields:
\begin{align*}
|C_2 \cap R(I(C_1))|    & \leq \frac{1}{r-1} \sum^{r}_{i=2}|C_i \cap R(I(C_1))| =\frac{1}{r-1} |\bigcup^{r}_{i=2} C_i \cap R(I(C_1))|\\
& = \frac{1}{r-1} |(X \smallsetminus C_1) \cap R(I(C_1))|  \leq \frac{r-2}{r-1}|I(C_1))|.
\end{align*}
and therefore
\begin{equation}
\label{equ:|C2|<|I(C1)|}
   |C_2 \cap R(I(C_1))| <|I(C_1)|.
\end{equation}

This, in turn, implies $|I(C_2)| \leq d$.

Next, we define the set $U$ to be plugged into Proposition $\ref{prop:C1..Cn}$.

For $i=1,2$ we divide $I(C_i)$ into two disjoint sets:
\[ S_i=\{j \in D: |C_i \cap R_j|=1\}\]
\[ M_i=\{j \in D: |C_i \cap R_j|>1\}\]
and get:
\begin{equation}
\label{equ:Ci>I(Ci)+Mi}
|C_i| \geq |I(C_i)|+|M_i|.
\end{equation}
Furthermore, for every subset $J$ of $D$:
\begin{equation} \label{eq:ineq}
|C_i \cap R(J)| \geq |I(C_i) \cap J|+|M_i \cap J|.
\end{equation}
Assume $|I(C_1)|=d-a\, (a \geq 0)$. From (\ref{equ:Ci>I(Ci)+Mi}) we obtain
\[ |M_1| \leq |C_1|-|I(C_1)| \leq d-(d-a)=a\]
and from ($\ref{equ:|C2|<|I(C1)|}$):
\[ |C_2 \cap R(I(C_1))| \leq d-a-1.\]
The set $S_1 \cap S_2$ can be divided into two disjoint subsets:
\[ U_1=\{ j \in S_1 \cap S_2: C_1 \text{ is lower than } C_2 \text{ on } R_j\},\]
\[ U_2=\{ j \in S_1 \cap S_2: C_2 \text{ is lower than } C_1 \text{ on } R_j\}.\]
If $|U_1| \geq |U_2|$ we define $U=U_1$, otherwise we define $U=U_2$.
In any case, $|U| \geq \frac{1}{2}|S_1 \cap S_2|$.

\begin{prop}
Under these notations
\[|I(C_1) \cap I(C_2) \smallsetminus U| \leq \lbrackk \frac{d-1}{2}\rbrackk .\]
\end{prop}
\begin{proof}
It suffices to show that
$2|I(C_1) \cap I(C_2) \smallsetminus U| \leq d-1$. Indeed:

\[
\begin{array}{ll}
2|I(C_1) \cap I(C_2) \smallsetminus U| & =  2|I(C_1) \cap I(C_2)|-2|U| \\
& \leq  2|I(C_1) \cap I(C_2)|-|S_1 \cap S_2| \\
& =     2|I(C_1) \cap I(C_2)|-|(I(C_1) \smallsetminus M_1) \cap S_2| \\
& =     2|I(C_1) \cap I(C_2)|-|I(C_1) \cap S_2|+|M_1 \cap S_2| \\
& \leq  2|I(C_1) \cap I(C_2)|-|I(C_1) \cap S_2|+|M_1| \\
&=|I(C_1) \cap I(C_2)|+(|I(C_1) \cap I(C_2)|-|I(C_1) \cap S_2|)+|M_1| \\
&=|I(C_1) \cap I(C_2)|+|I(C_1) \cap M_2|+|M_1|\\
& \leq_\text{by (\ref{eq:ineq})}  |C_2 \cap R(I(C_1))|+|M_1| \\
& \leq ( d-a-1)+a\\
& =  d-1.
\end{array}
\]
\end{proof}
To finish the proof of theorem $\ref{thm:mainThm}$, we choose for each index $i \in I(C_1)\cap I(C_2)\smallsetminus U$
a set $C_q\;(3 \leq q \leq r)$ that does not meet $R_i$, and call it $C(i)$. (Such a set does exist, since
$|(X \smallsetminus (C_1 \cup C_2)) \cap R_i|\leq r-3$.) Note that the sets $C(i)\;(i\in I(C_1)\cap I(C_2) \smallsetminus U)$
are not necessarily distinct.

Under these conditions, the sets $C_1,C_2,\{C(i): i \in I(C_1)\cap I(C_2) \smallsetminus U \}$ satisfy the assumptions
of Proposition $\ref{prop:C1..Cn}$ with $n \leq 2+ \lbrackk \frac{d-1}{2}\rbrackk = \lbrackk \frac{d+3}{2}\rbrackk $, and therefore $\conv C_1 \cap \conv C_2 \cap \bigcap_{i\in I(C_1)\cap I(C_2) \smallsetminus U} \conv C(i) = \emptyset$.

To sum it up, we have shown that for every $d \geq 2$ and $r \geq
\lbrackk\frac{d+3}{2}\rbrackk$ there is a set $X$ of $(d+1)(r-1)$ points in
$\mathbb{R}^d$, such that in any $r$-partition of $X$, there are
$\lbrackk\frac{d+3}{2}\rbrackk$ parts whose convex hulls have empty intersection.
This completes the proof of Theorem~\ref{thm:mainThm}.

\section{Proof of theorem \ref{thm:342example}: the case $d=3,r=4,k=2$}

Consider again the set $X$ described in the previous section for
$d=3$ and $r=4$ (three points on each of four rays emanating from the
origin in $\mathbb{R}^3$). The proof of Theorem \ref{thm:mainThm}
($T(3,4,3)=13$) shows that in every $4$-partition of $X$ there are
three parts whose convex hulls have empty intersection. It does
\textbf{not} show that there are two parts whose convex hulls do not
touch. In fact, $X$ does admit a partition into four triangles that
meet pairwise at the edges.

In this section we shall modify $X$ slightly by applying a suitable
small perturbation to two of its points. This will yield a set $X'
\subset \mathbb{R}^3$ of $12$ points, such that in every
$4$-partition of $X'$, two of the parts  have disjoint convex hulls,
and therefore $T(3,4,2)=13$.

Call a $4$-partition of $X$ (or of another set $X'$) "bad" if some
two of the parts have disjoint convex hulls. A partition is "good" if
the convex hulls of each two parts have a point in common.

\begin{prop}\label{rem:perturb}
    If $X=C_{0} \cup C_{1} \cup C_{2} \cup C_{3}$ is a "bad"
    4-partition of $X$, it remains "bad" if we apply a sufficiently small perturbation to
    $X$.
\end{prop}
\begin{proof}
    Assume, e.g., that $\conv C_0 \cap \conv C_1=\emptyset$. Choose a
    positive number $\delta$ such that:
    \begin{enumerate}
        \item $\| x-y \|> 2 \delta$ for all $x \in \conv C_{0}$, $y
            \in \conv C_{1}$.
        \item $\| a-b \|> 2 \delta$ for all $a,b \in X$, $a \neq
            b$.
    \end{enumerate}
    Replace each point $a \in X$ by a point $a'$ satisfying $\| a-a' \|<
    \delta$, and define $C'_{i} = \{a' : a \in C_{i} \}$ $(i=0,1,2,3)$.
    Since moving each point of $C_{0}$ (or $C_{1}$) by less than
    $\delta$ will cause every convex combination to change by less
    than $\delta$, we still have $\conv C'_{0} \cap \conv C'_{1}=\emptyset$

\end{proof}

Since $X$ has only finitely many partitions, there is a positive
$\delta$, such that moving each point of $X$ less than $\delta$ will
leave every "bad" partition of $X$ "bad".

With the help of the next three propositions we shall determine all
the "good" $4$-partitions of $X$. (There are exactly eighteen.) Then
we shall convert them all into "bad" partitions by slightly moving
just two points of $X$.

\begin{prop}
Let $\mathcal{C}=(C_{0} , C_{1} , C_{2} , C_{3})$  be a $4$-partition
of $X$. If $|I(C_{i})|<3$  for some $i$, then $\mathcal{C}$ is "bad".
\end{prop}
\begin{proof}
Assume, e.g., that $|I(C_0)|<3$.

If $|I(C_0)|=1$ then $C_0$ is included in a single ray, say $R_0$.
Since $|X \cap R_0|=3$, there is a part, say $C_1$, that misses $R_0$,
i.e., $I(C_1) \cap I(C_0)= \emptyset$. It  follows that $ \conv C_{0} \cap
 \conv C_{1}=\emptyset$ ,by Proposition~\ref{prop:supp_x}.

If $|I(C_0)|=2$, then $C_0$ is included in two rays, say $C_0 \subset
R_0 \cup R_1$. Consider three cases:
\begin{enumerate}
    \item If $|C_0| \geq 4$, then $|(X \smallsetminus C_0) \cap (R_0
        \cup R_1)| \leq 2$. Therefore some other part, say $C_1$,
        misses $R_0 \cup R_1$ entirely, and therefore $ \conv C_0 \cap
         \conv C_1= \emptyset$, again by
        Proposition~\ref{prop:supp_x}.
    \item If $|C_0|=3$, then $|(X \smallsetminus C_0) \cap (R_0 \cup
        R_1)|=3$. If some other part $C_i$ misses $R_0 \cup R_1$
        entirely, then $ \conv C_0 \cap  \conv C_i= \emptyset$, as above.
        If not, then  $|C_i \cap (R_0 \cup R_1)| =1$ for
        $i=1,2,3$. One of these two rays, say $R_0$, carries a
        single point of $C_0$ and one point each of two other
        parts, say $C_1$ and $C_2$. Applying
        Proposition~\ref{prop:C1..Cn} to $C_0$ and $C_1$, with
        $U= \{ 0 \} $, we find that $ \conv C_0 \cap  \conv C_1=
        \emptyset$.
    \item If $|C_0|=2$, then $|(X \smallsetminus C_0) \cap (R_0 \cup
        R_1)| =4$. If some other part misses $R_0 \cup R_1$,
        apply  Proposition~\ref{prop:supp_x} as above. If not,
        then $|C_i \cap (R_0 \cup R_1)| =1$ for some $i$, $1 \leq i
        \leq 3$. In this case, apply
        Proposition~\ref{prop:C1..Cn} to $C_0$ and $C_i$ as in
        the previous case.
\end{enumerate}
\end{proof}
Henceforth we assume that $|I(C_i)| \geq 3$ for $0 \leq i \leq 3$, and
therefore $|C_i| \geq 3$. Recalling that $|X|=12$ we conclude that $|C_i|=3$
for $0 \leq i \leq 3$. In other words, each part contains exactly $3$ points
from $3$ different rays.
It follows easily that no two parts occupy the same three rays, and therefore associating each part $C_i$ with the ray disjoint from $C_i$ is a bijection. For convenience, denote the ray disjoint from $C_i$ by $R_i$.

Now, let us look at two parts: $C_i,C_j$ where $I(C_i)= \{ j,k,l \}$ and $I(C_j)= \{ i,k,l \}$. ($\{ i,j, k,l \}=\{ 0,1, 2,3 \}$.)
\begin{prop}
\label{prop:CiLowerThanCj}
If $C_i$ is lower than $C_j$ on both $R_k$ and $R_l$, then $ \conv C_i \cap  \conv C_j=\emptyset$.
\end{prop}
\begin{proof}
This follows from Proposition~\ref{prop:C1..Cn} with $U= \{ k,l \}$.
\end{proof}

\begin{definition}
Denote the highest, middle and lowest points of $X$ on the ray $R_i$ by $p_{i,1},p_{i,2},p_{i,3}$ respectively.
\end{definition}

\begin{prop}
\label{prop:LowestMiddleHighest}
If $\mathcal{C}=(C_{0} , C_{1} , C_{2} , C_{3})$ is a "good" partition, then each set $C_i$ contains exactly one lowest point, one middle point, and one
highest point.
\end{prop}
\begin{proof}
If the lowest points of both $R_k$ and $R_l$ are in $C_i$, then, due to
Proposition~\ref{prop:CiLowerThanCj}, $\conv C_i \cap \conv C_j = \emptyset$, so $\mathcal{C}$ is  a "bad" partition. For
the same reason, $C_i$ doesn't contain two highest points. Therefore, the four lowest points are divided among the four parts, one for each, and the same is true for the highest points. This leaves no choice for the middle points, but to be divided evenly among the four parts.
\end{proof}

Let us denote by $H_{k,l}$ the plane $\aff(R_k \cup R_l)$ and by $H_{k,l}(j)$ the closed half-space bounded by $H_{k,l}$ that includes $R_j$. Then $H_{k,l}$ (weakly) separates the two triangles $\conv C_i$ and $\conv C_j$, with $\conv C_i \subset H_{k,l}(j)$ and $\conv C_j \subset H_{k,l}(i)$. Each of these two triangles meets $H_{k,l}$ in an edge, and these two edges cross iff $C_i$ is lower than $C_j$ in $R_k$ and higher than $C_j$ in $R_l$, or vice versa. We can (strictly) separate these two triangles by pushing a point of $C_i \cap H_{k,l}$ slightly into $ H_{k,l}(j)$, or a point of $C_j \cap H_{k,l}$ slightly into $ H_{k,l}(i)$.

\begin{figure}[h]
\centering
\includegraphics[width=0.8\textwidth]{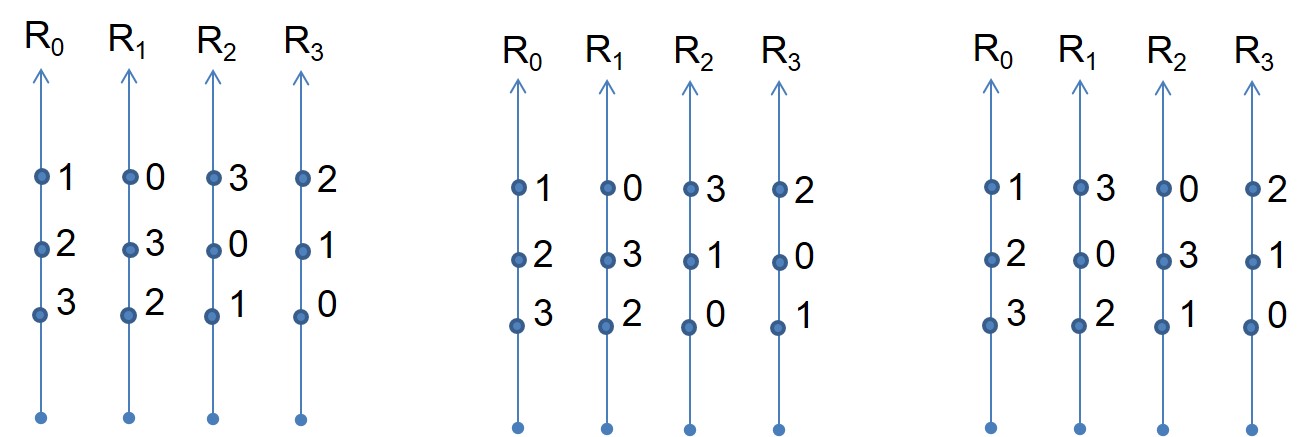}
\caption{Partitions}
\label{fig:partitionsOf4RaysInto4Parts}
\end{figure}

There are exactly three "good" partitions with $p_{0,1} \in C_1,\,p_{0,2} \in C_2$ and $p_{0,3} \in C_3$, as indicates in Figure \ref{fig:partitionsOf4RaysInto4Parts}.

The fifteen remaining "good" partitions are obtained by applying all permutations of $\{1,2,3\}$ (save the identity) to the three examples above. (Each permutation should act on all occurrences of the numbers $1,2,3$, including the subscripts in $R_1,R_2,R_3$.)

Now we are going to perturb the points of $X$ in order to avoid intersection in the "good" cases.

This will be a two-step process.  In step 1 we replace $p_{0,1}$
by another point $p'_{0,1}$, and call the resulting set $X'$. For
each $4$-partition  $\calc = (C_0, C_1, C_2, C_3)$ of $X$, we
denote by $\calc' = (C_0', C_1', C_2', C_3')$ the partition of
$X'$ obtained from $\calc$ when $p_{0,1}$  is replaced by
$p'_{0,1}$ in its part.  In step 2 we replace $p_{3,3}$ by another
point $p''_{3, 3}$, and call the resulting set $X''$. We denote by
$\calc'' = (C_0'', C_1'', C_2'', C_3'')$ the $4$-partition of
$X''$ obtained from $\calc'$ by substituting $p''_{3, 3}$ for
$p_{3, 3}$ in its part.

\begin{description}
\item [Step 1] Replace $p_{0,1}$ by $p'_{0,1}=p_{0, 1} +\delta' p_{3, 1},\;
\delta' > 0$.  Choose $\delta' $ sufficiently small, to ensure that
whenever $\calc$ is a "bad" $4$-partition of $X$, $\calc'$ will be a
"bad" $4$-partition of $X'$.  We claim that all "good"
partitions $\calc=(C_0, C_1, C_2, C_3)$ of $X$ (with $C_i \cap R_i =
\emptyset $ for $i = 0, 1, 2, 3)$ that satisfy $p_{0, 1} \in C_1$
or $p_{0, 1} \in C_2$ turn "bad".

Assume $p_{0, 1} \in C_1$. As mentioned above, $C_1$ and $C_3$ are (weakly) separated by $H_{0,2}$, with $C_1 \subset H_{0,2}(3),\,C_3 \subset H_{0,2}(1)$. Replacing  $p_{0, 1}$ by  $p'_{0, 1}$ amounts to pushing  $p_{0, 1}$ into $H_{0,2}(3)$, and thus strictly separating $C_1$ (now called $C_1'$) from $C_3$ (=$C_3'$). The same argument (with $C_1,H_{0,2},H_{0,2}(3),H_{0,2}(1)$ replaced by $C_2,H_{0,1},H_{0,1}(3),H_{0,1}(2)$ respectively) shows that if  $p_{0, 1} \in C_2$, then $\conv C_2' \cap \conv C_3'=\emptyset$.

%
%
%
%

We are still left with the "good" partitions $\calc = (C_0, C_1, C_2, C_3)$ of $X$ with $p_{0,1} \in C_3$. What can we say about $p_{3,3}$ in this case?
If $p_{0,2} \in C_1$, then $C_1$ is higher than $C_2$ on $R_0$, therefore
the highest point of $R_3$ cannot be in $C_1$. By
Proposition~\ref{prop:LowestMiddleHighest}, the
middle point of $R_3$ cannot be in $C_1$ (since the middle point
of $R_0$ is already in $C_1$). Therefore $p_{3,3} \in C_1$.
Similarly, if $p_{0,2} \in C_2$ then $p_{3,3} \in C_2$.

Conclusion:
$p_{3,3}$ is either in $C_1$ or in $C_2$.

Moreover, if  $p_{0, 1} \in C_3$, then  $p'_{0, 1} \in C_3'$, and therefore $C_0'=C_0$, $C_1'=C_1$ and $C_2'=C_2$. In other words, the parts $C_0,C_1$ and $C_2$ are unaffected by the perturbation in step 1.

\item [Step 2] Replace $p_{3,3}$ by $p''_{3,3}=p_{3,3}+ \delta'' p_{0,3}$, $ \delta''>0$.
(Again, $\delta''$ is small enough, so as to leave all "bad'' partitions of $X'$ "bad''.
Now we are going to show that if $p_{0,2}$ and $p_{3,3}$ are in $C_1$, then $\conv C_1'' \cap \conv C_0''=\emptyset$ (and similarly if $p_{0,2}$ and $p_{3,3}$ are in $C_2$, then $\conv C_2'' \cap \conv C_0''=\emptyset$).

Note that $C_1 \subset H_{2,3}(0),\,C_0 \subset H_{2,3}(1)$. If $p_{3,3}$ is in $C_1$, then the replacement $p_{3,3} \rightarrow p_{3,3}''$ anounts to pushing $p_{3,3}$ into $\inter H_{2,3}(0)$, and thus strictly separating $\conv C_1''$ from $\conv C_0''(=\conv C_0)$.

The same argument, with the roles of "$1$" and "$2$" interchaged, applies when $p_{3,3}$ is in $C_2$. 
%

\end{description}

After steps 1 and 2 we obtain a set $X''$ of 12 points in $\mathbb{R}^3$, such that for any partition of $X''$ into 4 disjoint parts, some two parts have disjoint
convex hulls. This completes the proof of Theorem~\ref{thm:342example}.

\section{Proof of theorem \ref{thm:secondThm}}

For this proof we use the same counter-example as in Theorem
$\ref{thm:mainThm}$, with an additional restriction. Recall that we started with
a simplex, centered at the origin, with $d+1$ vertices $p_0,\ldots,p_d$. For each vertex $p_i$ we defined $R_i$ to be the open ray
emanating from $\underline{0}$ through $p_i$. On each ray we chose
$r-1$ points. The set of all these chosen points is denoted by $X$.
$|X|=(d+1)(r-1)=T(d,r)-1$. The additional restriction in our case is that the
$r-1$ points on each ray $R_i$ ($i=0,\ldots,d$) be in "general
position", as detailed in the next paragraph.

For a subset $M \subset D$ ($D=\{0,1,\ldots,d\}$), define an $(M,X)$-selection $S$ to be a subset of $X$, of size $|M|$, consisting of exactly one point on each ray $R_j,\,j \in M$. The set  $X (\subset \cup \{R_i : i \in D\})$ is in "general position" if for any set $M \subset D$, $2 \leq |M| =m \leq d$, and for every $\bar{m}$ pairwise disjoint $(M,X)$-selections $S_1,\ldots,S_{\bar{m}}$, the intersection $\cap_{i=1}^{\bar{m}} \aff S_i$ is a single point if $\bar{m}=m$, and is empty if $\bar{m}>m$. (Since the maximum possible number of pairwise disjoint $(M,X)$-selections is just $r-1$, this condition applies only to sets $M \subset D$ of size $2 \leq |M| \leq \min \{r-1,d\}$.) A necessary and sufficient condition for this to happen is that if $S_i = \{ \lambda_{i,j}p_j : j\in M\}, \quad i=1,\ldots,\bar{m}$ and  $M=\{j_1,\ldots,j_m\}$, then
$\det \left(
                                   \begin{array}{ccc}
                                     \lambda_{1,j_1}^{-1} & \cdots & \lambda_{1,j_m}^{-1}  \\
                                     \vdots &   & \vdots \\
                                     \lambda_{m,j_1}^{-1} & \cdots & \lambda_{m,j_m}^{-1} \\
                                   \end{array}
                                 \right)\neq 0
$ if $\bar{m}=m$, and
 $\det \left(
                                   \begin{array}{cccc}
                                     \lambda_{1,j_1}^{-1} & \cdots & \lambda_{1,j_m}^{-1} & 1 \\
                                     \vdots &  & \vdots & \vdots \\
                                     \lambda_{m+1,j_1}^{-1} & \cdots & \lambda_{m+1,j_m}^{-1} & 1 \\
                                   \end{array}
                                 \right)\neq 0
$ if $\bar{m}=m+1$.

\[ \]
 We will show that if
$d<\frac{rk}{r-k}-1$, then for any $r$-partition $(C_1,\ldots,C_r)$ of $X$, some $k$ of the convex hulls $\conv C_1,\ldots,\conv C_r$ have empty intersection. 

We proceed by induction on $r$.  This will enable us to focus on
partitions $(C_1,\dots, C_r)$ of $X$ where each part $C_j$
misses at least one ray $R_i$.

For $r =2$ there is nothing to prove.  Now assume $r > 2$, and
suppose the theorem holds for $r-1$.  Let $(C_1,\dots, C_r)$ be an
$r$-partition of the set $X$ defined above.  If one of the parts,
say $C_r$, contains a point from each ray $R_i$, then we turn to
the induction hypothesis.  We delete $C_r$, define $\tilde X = X \smallsetminus C_r$, and consider the $(r-1)$-partition $(C_1,\dots, C_{r-1})$ of $\tilde X$.  Note that $\tilde X$ contains at most $r-2$ points on
each ray $R_i$ and is in "general position", like $X$.

If $2\le k < r-1$, apply the induction hypothesis:  By assumption,
$d < \frac{rk}{r-k} - 1$ and since $\frac{rk}{r-k} -1 <
\frac{(r-1)k}{r-1-k} - 1, \; \tilde X$ satisfies  the conditions
of the theorem, and therefore some $k$ of the convex hulls $ \conv C_j
(j = 1, \dots, r-1)$ have empty intersection.

If $k = r-1$, then $\bigcap^{r-1}_{j=1}  \conv C_j= \emptyset$.  Indeed
if $x \in\bigcap^{r-1}_{j=1}  \conv C_j$, then $\supp x \subset
\bigcap^{r-1}_{j=1} I(C_j)$ by Proposition $\ref{prop:supp_x}$.  But
$\bigcap^{r-1}_{j=1} I(C_j) = \emptyset$, since each ray $R_i$ is
missed by at least one of the parts $C_1,\dots, C_{r-1}$.  The
only point $x\in\bbr^d$ with $\supp x = \emptyset$ is the origin
$\underbar{0} $, but $\underbar{0} \notin  \conv C_j$ unless $I(C_j) =
D$.

From now on we assume that for every $j$, $I(C_j)\subsetneq D$.

We now prove the theorem. To do this we define a (weight)
function: given $k$ distinct parts (say $C_{j_1},\ldots,C_{j_k}$)
and  a ray $R_i$,  define: \begin{multline*}
W((C_{j_1},\ldots,C_{j_k}),R_i)  := \left\{\begin{array}{l}
0\,\,\,\,\,\,\,\,\,\,\,\,\, \text{if }\ R_i \cap C_{j_s}=\emptyset \text{ for some } s \in \{1,\ldots,k\}\\
\\
1+\#\{s:|C_{j_s}\cap R_i|>1\}\,\,\,\,\,\,\,\,\,\,\,\,\, \text{otherwise}\, .
\end{array}\right.\end{multline*}

In Section $\ref{subsec:Lower_Bound}$
we will show that if $\cap_{s=1}^k \conv C_{j_s}\neq \emptyset$ then
\begin{equation}\label{eq:buba} \sum_{i=0}^d W((C_{j_1},\ldots,C_{j_k}),R_i)\geq k.
\end{equation}

In Section
$\ref{subsec:Upper_Bound}$
we will show that for each $i\in D$:

\begin{equation}\label{eq:baba}\sum_{1\leq j_1<j_2<\ldots<j_k\leq r}W((C_{j_1},\ldots,C_{j_k}),R_i)\leq \binom{r-1}{k}.
\end{equation}

We use these two results to establish Theorem
$\ref{thm:secondThm}$.  If $\cap_{s=1}^k \conv C_{j_s}\neq \emptyset$ for all
$1\leq j_1<j_2<\ldots<j_k\leq r$, then from the inequalities
($\ref{eq:buba}$) and ($\ref{eq:baba}$) we conclude:

\begin{equation}\label{eq4.0.6}  k\binom{r}{k}\leq \sum_{1\leq j_1<j_2<\ldots<j_k\leq
r}\sum_{i=0}^d W((C_{j_1},\cdots,C_{j_k}),R_i)\leq
(d+1)\binom{r-1}{k}.\end{equation}

We thus obtain:

\[k\binom{r}{k}\leq (d+1)\binom{r-1}{k} \]

which is equivalent to $d\geq \frac{rk}{r-k}-1$, and the theorem
follows.

\subsection{A lower bound for the weight function $W$.}
\label{subsec:Lower_Bound}

Let $\{ C_j\}_{j\in J} \; ( J \subset\{ 1, \dots, r\}, |J| = k)$  be
a collection of $k$ parts.  We aim to show that if $\bigcap\{
 \conv C_j: j \in J\} \neq \emptyset$ then $\suml^d_{i = 0} W(\{
C_j\}_{j \in J}, R_i\} \ge k$.

For the weight $W(\{ C_j\}_{j \in J}, R_i)$ to be positive, each
of the parts $C_j \,(j \in J)$ must meet the ray $R_i$.  We say that
$R_i$ is a {\bf common} ray (for the given collection) if $R_i
\cap C_j \neq \emptyset$ for all $j \in J$.  For convenience we
define $I(J)= \bigcap_{j\in J} I(C_j)$ to be the set of indices of
the common rays.  Then the union of the common rays is just
$R(I(J))$.  Proposition $\ref{prop:index}$ below  says that the intersection of
the convex hulls $\bigcap \{  \conv C_j: j \in J\}$ depends only on the
intersections of the parts $C_j( j\in J)$ with the common rays.

\begin{prop}
\label{prop:index} For $J\subset \{1,\ldots,r\}$, if
$I(C_j)\subsetneq D$ for every $j\in J$ then
\[ \cap_{j\in J} \conv C_j=\cap_{j\in J}\conv(C_j\cap R(I(J))).\]
\end{prop}

\begin{proof}
The r.h.s. is clearly a subset of the l.h.s.. We
show that the l.h.s. is included in the r.h.s. as follows: suppose
$x\in \cap_{j\in J}\conv C_j$. Then, by Proposition
\ref{prop:supp_x}, $\supp x\subset I(J)$. By Proposition
\ref{prop:CcapR(supp_x)} it follows that $x\in \conv(C_j\cap
R(I(J)))$ for all $j\in J$.
\end{proof}

\begin{prop}
\label{prop:mutualrays}
 Given a set of $k$ parts,
$
\{C_j\}_{j \in J} (J \subset \{ 1, \dots, r\} , |J| = k)$, if $|I(J)| = m < k,$ and if each of the common rays $R_i \,(i \in
I(J))$ contains exactly one point on each of the $C_j$-s, then
$\bigcap_{j\in J} \conv C_j = \emptyset$.
\end{prop}
\begin{proof} Since $|I(J)| = m$, we have $m$ common rays spanning
an $m$-dimensional linear space.  Each of the sets $C_j \cap
R(I(J))$ consists of $m$ linearly independent points and therefore
spans a hyperplane in that space.  By the definition of "general
position" (see above), we have:
\[\bigcap_{j\in J} \aff (C_j\bigcap R(I(J))) = \emptyset.\]
\end{proof}

The following is a natural generalization of the last proposition:

\begin{prop}
\label{prop:m<k-t}
Given a set of $k$ parts, $\{C_j\}_{j\in J}$ ($J\subset
\{1,\ldots,r\}$, $|J|=k$), suppose $|I(J)|=m$ and denote by $t$
the number of parts among the $C_j$-s that contain more than one
point of at least one of the common rays. In this case, $m<k-t$
implies $ \cap_{j\in J} \conv C_j=\emptyset. $
\end{prop}
\begin{proof}
Divide $J$ into two subsets $S, T$ as follows:

$j\in S$ if $C_j$ meets each ray $R_i \,(i \in I(J))$ in a single
point.

$j\in T$ if $C_j$ meets at least one ray $R_i \,(i \in I(J))$ in
more than one point.

Then $|T| = t, |S| = k-t.$

By Proposition $\ref{prop:index}$,
\begin{align*}
\bigcap_{j\in J}  \conv C_j &= \bigcap_{j\in J}  \conv (C_j\cap R(I(J)))
\\
& \subseteq \bigcap_{j\in S}  \conv (C_j\cap R(I(J)))\\
& \subseteq \bigcap_{j\in S} \aff (C_j \cap R(I(J))).
\end{align*}
The last expression is the intersection of $k-t (> m)$ hyperplanes
in the $m$-dimensional space spanned by $R(I(J))$, which is empty
due to the "general position" of $X$. \end{proof}

Proposition $\ref{prop:m<k-t}$ implies inequality ($\ref{eq:buba}$). Indeed, from Proposition $\ref{prop:m<k-t}$ it follows that if $J=\{j_1,\ldots,j_k\}$  and $\cap_{j \in J}\conv C_j \neq \emptyset$, then $m+t \geq k$. But the l.h.s. of  ($\ref{eq:buba}$) is just $m+\# \{(i,s):i \in I(J),\,1 \leq s \leq k \text{ and } |C_{j_s} \cap R_i|>1\}$, which is $\geq m+t$.

\subsection{An upper bound for the weight function $W$.}
\label{subsec:Upper_Bound}
The weight of a ray $R_i$ is defined as:
\[
W(R_i):= \suml_{1\le j_1 < j_2 < \cdots < j_k\le r} W((C_{j_{1,
\dots,}} C_{j_k}), R_i).\]
 We will show that $W(R_i)$ is maximal when each point in $R_i$
 belongs to a different part $C_j$, i.e., when $|C_j\bigcap R_i| \le
1 $ for all $j$.  In that case it is clear that $W(R_i)  =
 {r-1 \choose k}$.

 Now assume $|C_j\cap R_i| > 1$ for some $j$, say $|C_1\cap R_i| >
 1$. Since $|X\cap R_i| = r -1$, there is another part, say
 $C_2$, that does not meet $R_i$ at all.  Choose one point $ x \in
 C_1\cap R_i$, and change the given partition $\mathcal{C} = (C_1, \dots,
 C_r)$ into $\mathcal{C'} = (C'_1,\dots, C'_r)$ as follows:
\[C'_1 = C_1 \smallsetminus \{ x \},\; C_2' = C_2\cup\{x\},\; C_j' = C_j
\;\; \text{for } \, 3 \le j \le r.\]

This change will increase the value of $W(R_i)$, or
leave it unaffected.  In fact, if $|C_1\cap R_i| > 2$, then\
\[W(\{ C'_j: j \in J\}, R_i) \ge W((C_j: j \in J\}, R_i)\]
for all $k$-subsets $J \subset D$.  If $|C_1\cap R_i| = 2$, define
$P_i = \{ j \in \{ 1, \dots, r\}: \: C_j\cap R_i \neq \emptyset\}$
and note that \[ W(\{C'_j : j \in J\}, R_i) = W(\{C_j: j \in J\},
R_i) - 1\] iff $J \subseteq P_i \;(|J| = k)$ and $1\in J$.  This
happens exactly ${|P_i|-1\choose k-1}$ times. On the other hand,

\[ W(\{ C'_j: j \in J\}, R_i) \ge  W(\{ C_j:j\in J\}, R_i)+1 \;
\; (=1) \]
 iff $2 \in J\;(|J|=k)$ and $J \smallsetminus\{ 2\} \subset P_i$.
  This happens exactly
${|P_i|\choose k-1}$ times. For all other $k$-sets $J\subset D$
there is no change at all. Since ${ |P_i|\choose k-1 } -
{|P_i|-1\choose  k-1} = {|P_i|-1 \choose k-2} \geq 0$, the total
change in $W(R_i)$ is nonnegative.

We can repeat this operation until all $r-1$  points of $X \cap
R_i$ belong to different parts $C_j$, in which case $W(R_i) =
{r-1\choose k}$.  Thus initially $W(R_i) \le {r-1\choose k}$, as
claimed in ($\ref{eq:baba}$).

If initially $|P_i| < r - 1$, then in the last step of the process
described above $|P_i|$ increases from $r-2$ to $r-1$.  In that
step $W(R_i)$ increases by ${|P_i |-1\choose k - 2} =
{r-3\choose k -2}$, which is strictly positive, since $0 \le k -
2 \le r-3$.  This shows that $W(R_i) = {r-1\choose k}$ iff the
$r-1$ points of $X\cap R_i$ belong to $r-1$ different parts $C_j$.

\begin{rmk}
\label{rmk:d+1=rk/(r-k)}
In case $d+1 = \frac{rk}{r-k}$ (or, equivalently, $k{r\choose k} =
(d+1) {r-1\choose k}$) we can repeat the arguments of the proof of
Theorem $\ref{thm:secondThm}$ and find that if $\calc = (C_1,\dots, C_r)$ is an
$r$-partition of $X$, and each $k$ of the convex hulls $\conv C_j$
have a point in common, then inequality  ($\ref{eq4.0.6}$) holds. (This is true if we assume that no part $C_j$ visits all $d+1$ rays $R_0,\ldots,R_d$. But if one part visits all rays, some $k$ of the convex hulls of the remaning $r-1$ parts have empty intersection. This is shown in detail in the earlier part of the proof of  Theorem $\ref{thm:secondThm}$.) Since
$k{r\choose k} = (d+1) {r-1\choose k}$, both inequalities in
($\ref{eq4.0.6}$) must hold as equalities.  In view of ($\ref{eq:baba}$), this
implies that $W(R_i) = \suml_{1\le j_1 < j_2 < \cdots < j_k \le r}
W((C_{j_1} \dots, C_{j_k}),R_i) = {r-1\choose k}$ for $i = 0, 1, \dots,
d$.  This, in turn, implies that each ray $R_i$ carries $r-1$
points of $X$ that belong to $r-1$ different parts.  One can
easily deduce that for each $k$ distinct parts $C_{j_1}, \dots,
C_{j_k}$, there are exactly $k$ rays $R_i$ that intersect each of
these parts, and therefore the convex hulls $\conv C_{j_1},\dots,
\conv C_{j_k}$ intersect in a single point.  In these cases there is
some hope to transform $X$ by a small perturbation into a "bad"
set $X'$, such that in any $r$-partition of $X'$ there are some
$k$ parts whose convex hulls have empty intersection.  In the next
section we shall do this in the case $d=5, r=3, k=2$.

\end{rmk}


\section{Proof of Theorem \ref{thm:532example}: the case $d=5$, $r=3$, $k=2$}

We start with the usual construction of $X$ (as in Theorem \ref{thm:mainThm}): six rays $R_i,\ i=0,\ldots,5$ in $\fR^5$, with two points chosen on each ray. The set $X$ contains $12=T(5,3)-1$ points. As before, we denote by $p_{i1}$ the upper point and by by $p_{i2}$ the lower point of $X$ on $R_i$.

Consider partitions of $X$ into three parts $X=C_1\cup C_2\cup C_3$. As in the case $d=3,r=4,k=2$, we say that a partition is "bad" if there are two parts $C_i,C_j$ such that $\conv C_i \cap \conv C_j=\emptyset$, and is "good" if for every $1\leq i<j\leq 3$, $\conv C_i \cap \conv C_j \neq\emptyset$.
Our aim in this section is to show that by moving three points of $X$ we can turn all "good" partitions into "bad" ones.

\begin{prop}
\label{prop:goodPartitions532}
A partition ($C_1,C_2,C_3$) is "good" iff for every $1\leq i<j\leq 3$ there exists exactly one ray where the lower point belongs to $C_j$ and the higher point belongs to $C_i$ and exactly one ray where the lower point belongs to $C_i$ and the higher point belongs to $C_j$.
\end{prop}

\begin{proof}

Let $(C_1,C_2,C_3)$ be a "good" partition. By Remark \ref{rmk:d+1=rk/(r-k)} above (with $d=5,\,r=3,\,k=2$), we find that for each two distinct parts, $C_i$ and $C_j$, there are exactly two rays that meet both $C_i$ and $C_j$. From Proposition \ref{prop:C1..Cn} with $n=2$ we conclude that $C_i$ must be lower than $C_j$ on one of these rays, and higher  than $C_j$ on the other.
\end{proof}

We see that the "good" partitions are those where for every $1\leq i\leq 3$, $|C_i|=|I(C_i)|=4$, and for every $1\leq i<j\leq 3$, $|I(C_i)\cap I(C_j)|=2$. Indeed, in these partitions $\conv C_i\cap \conv C_j\neq \emptyset$. Let us draw the following example of a "good" partition:
\begin{center}
\begin{picture}(100,70)(0,0)
\put(0,0){\usebox{\ry}} \put(-3,45){$R_0$}
\put(20,0){\usebox{\yr}} \put(17,45){$R_1$}
\put(40,0){\usebox{\yb}} \put(37,45){$R_2$}
\put(60,0){\usebox{\by}} \put(57,45){$R_3$}
\put(80,0){\usebox{\br}} \put(77,45){$R_4$}
\put(100,0){\usebox{\rb}} \put(97,45){$R_5$}
\end{picture}
\end{center}
Each ray is represented by a line starting at the origin $\underline{0}$ (the six different origins of the rays in the picture should of course be identified as one), and we denote the different parts of the partition by colors: yellow, blue and red. Proposition \ref{prop:goodPartitions532} implies that \emph{any} "good" partition is the same as this partition up to a permutation of the rays.
In fact, given any (other) "good" partition $(C_1, C_2, C_3)$, we can assign colors to the three parts (uniquely) in such a way that
$p_{0, 1}$ (the higher point on $R_0$) becomes red, $p_{0, 2}$ (the
lower point on $R_0$) becomes yellow, and the third part becomes blue.
Under this coloring, every "good" partition is obtained from the
example shown above by a permutation of $R_1,\ldots, R_5$. Thus
there are exactly $120(=5!)$ "good" partitions.
We shall dedicate some effort to study the example shown above.

Before continuing with the proof, recall the non-negative representation of a point (equation (\ref{equ:positRep})): if $p_0,\ldots,p_5$ are the vertices of a simplex whose center is at the origin ($\sum_i p_i =\underline{0}$), then
each point $x \in \bbr^5$ can be uniquely represented as $x= \sum \xi_i p_i$, where $\min (\xi_i)=0$.  A hyperplane in $\fR^5$ will be represented by an equation
\[ H=\{(\xi_0,\ldots,\xi_5)|\sum_{i=0}^5 a_i \xi_i=\alpha \}\]
where we demand that $\sum_i a_i =0$, so that the hyperplane $H$ be well defined. Indeed, if $(\xi_0,\ldots,\xi_5)=(\xi_0+\lambda,\ldots,\xi_5+\lambda)$ are two different representations of a point $x\in\fR^5$, then $ \sum_i a_i \xi_i=\sum_i a_i(\xi_i+\lambda) $ iff $\sum_i a_i=0$.

Returning to our example of a "good" partition, we notice that the convex hull of each of the color sets is a three dimensional simplex (inside $\fR^5$), and that any two of these simplices touch at a single point that lies in the relative interior of an edge of each of them. For example, the point of intersection of the red symplex and the yellow simplex is in $\text{span}\{R_0,R_1\}$ as depicted below.

\begin{figure}[h]
    \centering

    \begin{center}
\begin{picture}(100,70)(-50,0)
\put(0,0){\circle{5}}
\put(0,0){\line(1,1){60}}
\put(0,0){\line(-1,1){60}}
\put(-55,60){$R_0$}
\put(47,60){$R_1$}
\put(30,30){\line(-6,1){70}}
\put(37,30){red}
\put(-58,43){red}
\put(50,50){\line(-6,-2){70}}
\put(53,47){yellow}
\put(-50,20){yellow}
\end{picture}
\caption{The intersection of the red simplex and yellow simplex\label{fig:RedYelloyIntersectInAffR_0R_1}}
\end{center}

\end{figure}
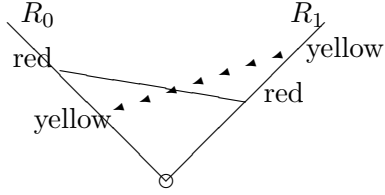

Moreover, any hyperplane of the form:
\[ H=\{(\xi_0,\ldots,\xi_5)|\quad a\xi_2+b\xi_3=c\xi_4+d\xi_5 \}\]
where $a,b,c,d>0 $ and $a+b=c+d$, (weakly) separates the yellow simplex from the red simplex: the red simplex lies in $H\cup H_+$ and the yellow simplex lies in $H\cup H_-$ ($\text{span}\{R_0,R_1\}$ is included in $H$). Here $H_{+},H_{-}$ are the open half spaces defined by $H_+=\{(\xi_0,\ldots,\xi_5)|\quad c\xi_2+d\xi_3<a\xi_4+b\xi_5 \}$ and $H_-=\{(\xi_0,\ldots,\xi_5)|\quad c\xi_2+d\xi_3>a\xi_4+b\xi_5 \}$.

We fix such a hyperplane $H$ by, say, choosing $a = b = c = d =
1$.  Then we apply a small perturbation to the point $p_{0, 1} \in
X$, trying to separate the  red simplex from the yellow one.  To
be concrete, we define $p'_{0, 1} = p_{0, 1} + \varepsilon\vec u$,
where $\vec u= (0, 1, 2, 3, 4, 5)$, and $\varepsilon$ is a
positive number, sufficiently small, so as to prevent the "bad"
partitions from becoming "good" (see Proposition $\ref{rem:perturb}$).

\begin{figure}[h]
    \centering

\begin{center}
\begin{picture}(100,70)(-50,0)
\put(0,0){\circle{5}}
\put(0,0){\line(1,1){60}}
\put(0,0){\line(-1,1){60}}
\put(-55,60){$R_0$}
\put(47,60){$R_1$}

\put(37,30){red}
\put(50,50){\line(-6,-2){70}}
\put(50,50){\circle{3}}
\put(-25,25){\circle{3}}
\put(30,30){\circle{3}}
\put(53,47){yellow}
\put(-55,20){yellow}
\end{picture}
\end{center}
\caption{The affine hull of $R_0 \cup R_1$ after the change of $p_{0,1}$\label{fig:RedAndYelloyNotIntersect}}
\end{figure}
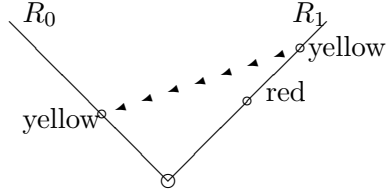

Now $p'_{0, 1} \in H_+$, and the red simplex lies in $H\cup H_+$,
with only one vertex in $H$. The yellow simplex lies, as before,
in $H\cup H_-$, with only one edge in $H$.  The red
vertex (in $H$) does not lie on the yellow edge (in $H$), and
therefore the yellow simplex and the red simplex are now disjoint. (See Fig. \ref{fig:RedAndYelloyNotIntersect}.)

The following two propositions yield sufficient conditions for any
perturbation to separate the convex hulls of two color sets.  They
will show that the perturbation described above separates the red
simplex from the yellow simplex in 100 out of the 120 "good"
partitions under consideration.

In order to spoil the remaining 20 "good" partitions we shall
need two more perturbations.

Let $\mathcal{C}$ be a partition of $X$ into three parts $C_1, C_2, C_3$
that conform to the description in Proposition $\ref{prop:goodPartitions532}$.

Assume
\begin{equation}
\label{equ:RaysOfC1_C2_C3}
\begin{cases} I(C_1) \cap I(C_3) = \{ i, j\}\\
I(C_2) \cap I(C_3) = \{ k, \ell\}\\
I(C_1) \cap I(C_2) = \{ m , n\}\end{cases}
\end{equation}
\[(\{ i, j, k , l , m, n\} = \{ 0, 1, 2, 3, 4, 5\}).\]
\begin{prop}
\label{prop:A_separate_hyperline}
Let $a, b, c, d$ be positive numbers, $a+b = c+d$.
Define:
\[ H_+ = \{ (\xi_0, \ldots, \xi_5): a \xi_i + b\xi_j > c\xi_k +d \xi_\ell\}\]

\[ H = \{ (\xi_0, \ldots, \xi_5): a \xi_i + b\xi_j = c\xi_k + d \xi_\ell\}\]

\[ H_- = \{( \xi_0, \ldots, \xi_5): a \xi_i + b\xi_j < c\xi_k + d
\xi_\ell\}.\] Then the hyperplane $H$ (weakly) separates $C_1$
from $C_2$, with $C_1 \subset H_+ \cup H, C_2\subset H_- \cup H$
and $C_1\cap H = C_1 \cap (R_m \cup R_n)$, $C_2\cap H = C_2 \cap
(R_m \cup R_n)$.
\end{prop}
\begin{proof} $R_i \cup R_j \subset \{ \underline{0}\} \cup H_+, \; R_k \cup
R_\ell \subset \{\underline{0}\} \cup H_-$, \ $R_m \cup R_n \subset H$.
\end{proof}
The following proposition deals with a perturbation of a point $p$ by defining $p'= p + \varepsilon \vec u$ ($\vec{u}=(u_0,u_1,u_2,u_3,u_4,u_5)$).
\begin{prop}
\label{prop:U_that_separates_c1_from_c2}
If $\max \{ u_i, u_j\} > \min \{ u_k, u_l\}$ (see ($\ref{equ:RaysOfC1_C2_C3}$)), and $
p\in C_1\cap (R_m\cup R_n)$, then the perturbation
$p\rightsquigarrow p'= p + \varepsilon \vec u \;(\varepsilon > 0 \text{ small})$ separates
$\conv C_1$ from $\conv C_2$.
\end{prop}
\begin{proof}
Assume, w.l.o.g., that $u_i\leq u_j$ and $u_k \leq u_l$. Define
\[
 \delta = \min \{ |u_\mu - u_\nu|: 0 \leq \mu,\nu \leq5,\,u_\mu \neq u_\nu\} \]
\[ \Delta = \max \{ u_\mu - u_\nu: 0 \leq \mu,\nu \leq5,\,u_\mu \neq u_\nu\}
\]
The assumption that $\max \{u_i,u_j\}>\min \{u_k,u_l\}$ implies $0< \delta \leq \Delta$. Now $u_j-u_k=\max \{u_i,u_j\}-\min \{u_k,u_l\}\geq \delta$, whereas $u_i-u_l\geq -\Delta$. It follows that $2\Delta (u_j-u_k)+\delta (u_i-u_l) \geq 2\Delta \delta - \delta \Delta = \Delta \delta >0$.

Now apply Proposition $\ref {prop:A_separate_hyperline}$ with $a=d=\delta,\,b=c=2\Delta$
 to obtain a hyperplane $H$ that weakly
separates $C_1$ from $C_2$: $C_1 \subset H_+ \cup H, \ C_2 \subset
H_- \cup H$, and $\conv C_1 \cap \conv C_2 = \conv (C_1 \cap H) \cap
\conv (C_2 \cap H)$ is the crossing point of the two segments
 $\big[p_{m, 1}, p_{n,2}\big]$ and $ \big[ p_{n, 1}, p_{m, 2}\big]$.
 (An example of that in given in Fig. \ref{fig:RedYelloyIntersectInAffR_0R_1} above.)

Since $p \in H$ and $\vec u \in H_+ (\delta u_i + 2\Delta u_j > 2\Delta u_k + \delta {u_\ell})$, we find that $p' = p+\varepsilon \vec u \in
H_+$. If we define $C_1' = C_1\smallsetminus \{ p\} \cup \{ p'\}$,  then
the simplex $\conv C'_1 $ lies in $H\cup H_+$ with only one vertex in
$H$.  This vertex misses the segment $H\cap \conv C_2$, and
therefore $\conv C_1' \cap \conv C_2 = \emptyset$.  (See Fig. \ref{fig:RedAndYelloyNotIntersect}.)
\end{proof}

{\bf First Step:} Applying Proposition $\ref{prop:U_that_separates_c1_from_c2}$ with $\vec u = (0, 1, 2, 3, 4, 5)$ (i.e.,
$u_\nu = \nu$ for all $\nu$), $C_1$ = red, $C_2$ = yellow, $m=0$ and
$p = p_{0, 1}$, we find that the first perturbation described
above does indeed separate the red simplex from the yellow
one, unless $\max \{ i, j\} < \min \{ k, \ell\}$, i.e., unless the
two red-blue rays precede the two yellow-blue rays.  These exceptional
partitions are exactly the 20 partitions shown below. (Each
figure represents four partitions, which differ only in the
internal ordering of the red-blue rays and the yellow-blue rays.)

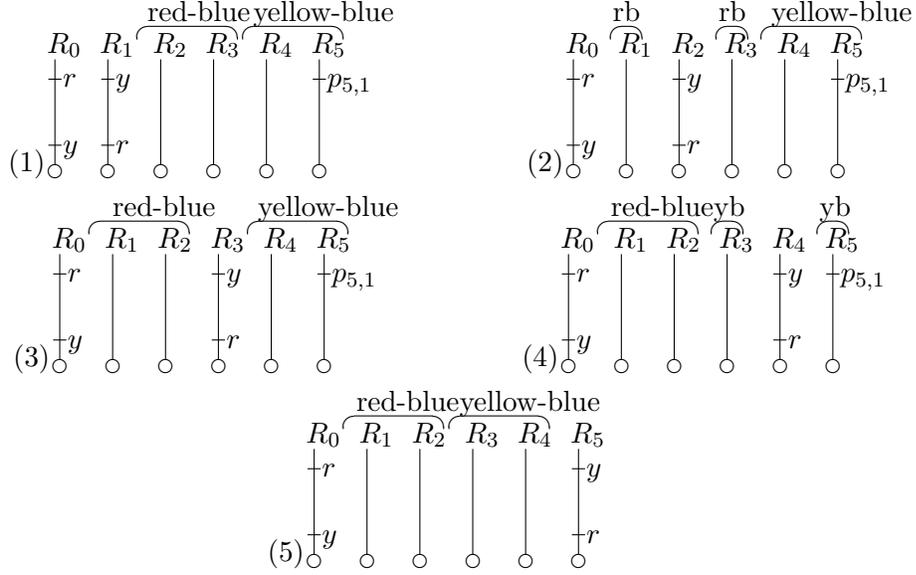
\begin{figure}[h]

(1)
\begin{picture}(100,70)(0,0)
\put(0,0){\usebox{\ry}} \put(-3,45){$R_0$}
\put(20,0){\usebox{\yr}} \put(17,45){$R_1$}
\put(40,0){\usebox{\simplebar}} \put(37,45){$R_2$}
\put(60,0){\usebox{\simplebar}} \put(57,45){$R_3$}
\put(35,57){red-blue}\put(50,50){\oval(38,10)[t]}
\put(80,0){\usebox{\simplebar}} \put(77,45){$R_4$}
\put(100,0){\usebox{\rightmost}} \put(97,45){$R_5$}
\put(75,57){yellow-blue}\put(90,50){\oval(38,10)[t]}
\end{picture} \hspace{2.5cm}
(2)
\begin{picture}(100,70)(0,0)
\put(0,0){\usebox{\ry}} \put(-3,45){$R_0$}
\put(40,0){\usebox{\yr}} \put(17,45){$R_1$}
\put(20,0){\usebox{\simplebar}} \put(37,45){$R_2$}
\put(60,0){\usebox{\simplebar}} \put(57,45){$R_3$}
\put(15,57){rb}\put(20,50){\oval(12,10)[t]}
\put(55,57){rb}\put(60,50){\oval(12,10)[t]}
\put(80,0){\usebox{\simplebar}} \put(77,45){$R_4$}
\put(100,0){\usebox{\rightmost}} \put(97,45){$R_5$}
\put(75,57){yellow-blue}\put(90,50){\oval(38,10)[t]}
\end{picture}\\
(3)
\begin{picture}(100,70)(0,0)
\put(0,0){\usebox{\ry}} \put(-3,45){$R_0$}
\put(60,0){\usebox{\yr}} \put(17,45){$R_1$}
\put(20,0){\usebox{\simplebar}} \put(37,45){$R_2$}
\put(40,0){\usebox{\simplebar}} \put(57,45){$R_3$}
\put(20,57){red-blue}\put(30,50){\oval(38,10)[t]}
\put(80,0){\usebox{\simplebar}} \put(77,45){$R_4$}
\put(100,0){\usebox{\rightmost}} \put(97,45){$R_5$}
\put(75,57){yellow-blue}\put(90,50){\oval(38,10)[t]}
\end{picture}\hspace{2.5cm}
(4)
\begin{picture}(100,70)(0,0)
\put(0,0){\usebox{\ry}} \put(-3,45){$R_0$}
\put(80,0){\usebox{\yr}} \put(17,45){$R_1$}
\put(20,0){\usebox{\simplebar}} \put(37,45){$R_2$}
\put(40,0){\usebox{\simplebar}} \put(57,45){$R_3$}
\put(16,57){red-blue}\put(30,50){\oval(38,10)[t]}
\put(95,57){yb}\put(100,50){\oval(12,10)[t]}
\put(55,57){yb}\put(60,50){\oval(12,10)[t]}
\put(60,0){\usebox{\simplebar}} \put(77,45){$R_4$}
\put(100,0){\usebox{\rightmost}} \put(97,45){$R_5$}
\end{picture}\\
(5)
\begin{picture}(100,70)(0,0)
\put(0,0){\usebox{\ry}} \put(-3,45){$R_0$}
\put(100,0){\usebox{\yr}} \put(17,45){$R_1$}
\put(20,0){\usebox{\simplebar}} \put(37,45){$R_2$}
\put(40,0){\usebox{\simplebar}} \put(57,45){$R_3$}
\put(16,57){red-blue}\put(30,50){\oval(38,10)[t]}
\put(80,0){\usebox{\simplebar}} \put(77,45){$R_4$}
 \put(97,45){$R_5$}
\put(60,0){\usebox{\simplebar}}
\put(55,57){yellow-blue}\put(70,50){\oval(38,10)[t]}
\end{picture}
   \caption{The remaining 20 partitions}
    \label{fig:TheRemaining20Partitions}
\end{figure}

\begin{rmk}
We are going to apply Proposition \ref{prop:U_that_separates_c1_from_c2}  three times in a row. After the first step, the red point $p'_{0,1}$ (that replaces $p_{0,1}$) does not lie any more exactly on the ray $R_0$. But this does not really matter, since in the second step we separate the blue simplex from the yellow one. After the second step, already two points ($p'_{0,1}$ and $p'_{5,1}$) do not lie exactly on the coresponding rays $R_0$ and $R_5$. But the third step, where we try to separate the red simplex from the blue one, uses a vector $u''(=(1,0,0,0,1,0))$ that is known in advance. This means that we can make a finite list of hyperplanes $H$ that might be used to (weakly) separate the red simplex from the blue one (before the third perturbation). (Actually four different hyperplanes will suffice. See proof of Prop \ref{prop:U_that_separates_c1_from_c2}). By choosing $\varepsilon$ and $\varepsilon'$ sufficiently small, we can make sure that $p_{0,1}$ and $p'_{0,1}$ (and, if $p_{5,1}$ is blue, also $p_{5,1}$ and $p'_{5,1}$) lie on the same side of each of those hyperplanes $H$.
\end{rmk}
{\bf Second Step:} Replace the point $p_{5, 1}$ on $R_5$ by $p'_{5,
1} = p_{5, 1} + \varepsilon' \vec u'$, where $\vec u' = (5,4, 3, 2, 1, 0)$
and $\varepsilon'$ is a sufficiently small positive number, so as to
leave all "bad" partitions  "bad".

Our aim now is to separate the blue simplex from the yellow
simplex.  This will certainly fail in the last four cases (5),
where $R_5$ is yellow-red, and the two yellow-blue rays ($R_3$ and
$R_4$) remain untouched.

Now apply Proposition $\ref{prop:U_that_separates_c1_from_c2}$ with $m=5$, with $\vec u'$ instead of $\vec
u, C_1$ = blue and $C_2$= yellow when $p_{5,1}$ is blue, or $C_1$
= yellow and $C_2$ = blue when $p_{5, 1}$ is yellow, to show that
this second perturbation does indeed separate the blue simplex
from the yellow one. This will succeed in 14 out of the remaining
16 cases ((1)-(4)). It will fail only on the two
subcases of (1), where $p_{5, 1}$ is colored blue (see Fig \ref{fig:TheRemaining20Partitions}),
since in these cases (only) the two red-yellow rays precede the two
red-blue rays.
\begin{figure}[h]

        \begin{enumerate}
\item Case (1):

\begin{center}
\begin{picture}(100,70)(0,0)
\put(0,0){\usebox{\ry}} \put(-3,45){$R_0$}
\put(20,0){\usebox{\yr}} \put(17,45){$R_1$}
\put(40,0){\usebox{\simplebar}} \put(37,45){$R_2$}
\put(60,0){\usebox{\simplebar}} \put(57,45){$R_3$}
\put(35,57){red-blue}\put(50,50){\oval(38,10)[t]}
\put(80,0){\usebox{\yb}} \put(77,45){$R_4$}
\put(100,0){\usebox{\by}} \put(97,45){$R_5$}
\put(75,57){yellow-blue}\put(90,50){\oval(38,10)[t]}
\end{picture}
\end{center}
\item Case (5):

\begin{center}
\begin{picture}(100,70)(0,0)
\put(0,0){\usebox{\ry}} \put(-3,45){$R_0$}
\put(100,0){\usebox{\yr}} \put(17,45){$R_1$}
\put(20,0){\usebox{\simplebar}} \put(37,45){$R_2$}
\put(40,0){\usebox{\simplebar}} \put(57,45){$R_3$}
\put(16,57){red-blue}\put(30,50){\oval(38,10)[t]}
\put(80,0){\usebox{\simplebar}} \put(77,45){$R_4$}
 \put(97,45){$R_5$}
\put(60,0){\usebox{\simplebar}}
\put(55,57){yellow-blue}\put(70,50){\oval(38,10)[t]}
\end{picture}
\end{center}
\end{enumerate}
    \caption{The 6 partitions that are still "good"}
    \label{fig:TheLast8Partitions}
\end{figure}

In the remaining cases (Fig. \ref{fig:TheLast8Partitions}) we shall try to separate the
blue simplex from the red simplex by means of a third
perturbation.  Please note that in these remaining cases the two
red-blue rays $(R_2$ and $R_3$ in (1), $R_1$ and $R_2$ in (2)) have
not been affected by the first two perturbations.

{\bf Third step:} Replace the point $p_{2, 1}$ on $R_2$ by $p_{2, 1}''
= p_{2, 1} + \varepsilon'' \vec u'', (\varepsilon'' > 0)$ where $\vec
u '' = (1, 0, 0, 0, 1, 0)$ and $\varepsilon ''$ is sufficiently small,
so as to leave all "bad" partitions "bad".

Apply again Proposition $\ref{prop:U_that_separates_c1_from_c2}$, with $m = 2, \vec u''$ instead of $\vec u$, $C_1$= blue and $C_2 $ = red when $p_{2, 1}$ is blue, or $C_1$
= red and $C_2$ = blue when $p_{2, 1}$ is red, to find that the
third perturbation does indeed separate the blue simplex from the
red one in all the six remaining cases.

This concludes the proof of Theorem \ref{thm:532example}: we began with a set $X$ of 12 points in $\fR^5$ and perturbed three of them to get a set $X''$ of 12 points in $\fR^5$, such that in any $3$-partition of $X''$ there are two parts whose convex hulls do not meet.

\section{Conclusion}This paper is devoted to the proof of parts of Reay's conjecture
($T(d,r,k)=T(d,r)$ for $2 \leq k \leq \min (d,r-1)$). The meaning of this conjecture (for specified values
of $d,r$ and $k$) is just this: there is a subset $X$ of $\bbr^d,\! |X|=T(d,r)-1\!(=(d+1)(r-1))$, such that in \textbf{every}
$r$-partition of $X$ ($X=C_1 \cup \cdots \cup C_r$) there are some $k$ parts whose convex hulls have empty intersection.
The conjecture is meaningful for all triples $(d,r,k)$ of values that satisfy $2 \leq k < d+1$ and $k<r$.
We prove the conjecture whenever $k+1 \leq d+1 \leq 2k-1$ (Theorem $\ref{thm:mainThm}$). When $2k \leq d+1 < k(k+1)$ we
prove the conjecture for $k<r<\frac{d+1}{d+1-k}k$ (see Remark following Theorem $\ref{thm:532example}$),
and also in the two special cases $(d,r,k)=(3,4,2)$ and $(d,r,k)=(5,3,2)$.
In all cases, the set $X$ is a variation, specialization or perturbation of the same example:
$d+1$ rays that emanate from the origin and positively span $\bbr^d$, with $r-1$ points chosen on each ray.

Unfortunately,we were unable to disprove Reay's conjecture for any admissible triple $(d,r,k)$.
It is conceivable, though, that $T(d,r,k)<T(d,r)$ holds for any given values of $r$ and $k$ ($2 \leq k <r$),
provided $d$ is large enough. In particular, one might try to show that $T(d,3,2)<T(d,3)=2d+3$ from some $d$
onward (maybe already for $d\geq 6$).

Note that claims concerning $T(d,r,2)$ ($k=2$) are actually statements about Radon partitions.
Radon partitions are much better understood an easier to handle than Tverberg $k$-partitions for $k\geq 3$.


\end{document}